\documentclass{notices}
\usepackage{amsfonts,amssymb,amsmath,amscd}
\usepackage{xcolor}
\usepackage{graphicx}
\usepackage{amssymb}
\usepackage{epstopdf}
\usepackage{xcolor}
\usepackage{ulem}
\usepackage{cancel}
\usepackage[percent]{overpic}
\usepackage{stackengine}
\usepackage{verbatim}

\title{Bernard Maskit Memorial Tribute\footnote{A shorter version of this paper appeared in the August 2025 issue of {\it The Notices}.\cite{Not1}}}

\author{James W. Anderson
  \affil{James W. Anderson is a professor of mathematics at 
 the University of Southampton UK. His email address is 
J.W.Anderson@soton.ac.uk}
\and
Ara Basmajian
  \affil{
    Ara Basmajian is a  professor of mathematics at the Graduate Center and Hunter College, CUNY. His  email address is abasmajian@gc.cuny.edu
    }
   \and
   Ruben Hidalgo
  \affil{Ruben Hidalgo is a  professor of mathematics  at 
  Universidad de La Frontera, Temuco. His email address is 
 ruben.hidalgo@ufrontera.cl}
  \and
  Perry Susskind
  \affil{Perry Susskind is  a professor of mathematics at 
  Connecticut College.  His  email address is pdsus@conncoll.edu
   }
    \and
  Edward C. Taylor
  \affil{Edward C. Taylor, a retired professor of mathematics, passed away on March 10, 2025}
}

\begin{document}

\maketitle

\section*{His life and work}

 Bernard Maskit, who we all called Bernie,  was born on May 27, 1935 and died on March 15, 2024 from exacerbations and complications of COVID-19.  He grew up in New York City, attended DeWitt Clinton High School in The Bronx, and did his undergraduate degree at New York University.  Bernie completed his doctoral thesis, \lq\lq On Klein's Combination Theorem,\rq\rq \ in 1964 at NYU under famed complex analyst  Lipman Bers (\cite{MR1306867}, \cite{MR1476998}).  Bernie's professional life began with a post-doc at the Institute for Advanced Study in 1964-65, followed by an assistant professorship at MIT from 1965 to 1972.  He then moved to S.U.N.Y. at Stony Brook (now Stony Brook University) where he remained until his retirement in 2008 as professor emeritus.  Bernie gave an invited talk about Kleinian groups at the 1974 International Congress of Mathematicians, and became an inaugural fellow of the American Mathematical Society in 2012.   He had twelve doctoral students and   45 descendants as of February 2025\footnote{https://genealogy.math.ndsu.nodak.edu}.
 
Though the entirety of his mathematical work stands on its own as serious and consequential, Maskit's mathematical oeuvre is also, importantly, an essential bridge on a long  path that stretches back to the classification of the topological types of surfaces, beginning with M\"{o}bius in 1861, along with the work of Klein, Fuchs, Koebe and Poincar\'{e}, in which the uniformization theorem was proved.  That is, any Riemann surface\footnote{A Riemann surface is a 1-dimensional complex manifold that may be thought of as a surface along with a conformal structure or a hyperbolic metric structure.}  $S$ has a regular holomorphic cover $\pi:D\rightarrow D/G=S$, where $D$ is one of the Riemann sphere $\hat{\mathbb{C}} =\mathbb{C} \cup \{\infty\}$, the complex plane $\mathbb{C}$, or the upper half-plane $\mathbb{H}^2$, and $G$ is a discrete group of M\"{o}bius transformations acting properly discontinuously\footnote{  The action of $G$ on   $\Omega(G)$ is {\it properly discontinuous} if any compact set 
$K \subset \Omega(G)$ is moved disjointly away from itself for all but finitely many elements of $G$.   } on $D$.  The other end of this long path is the geometrization conjecture of William P. Thurston, the proof of which was completed by Grigori Perelman in 2003.  The latter is a result of celebrated significance in geometry and topology that occurred at the end of the 20th century and beginning of the 21st, and has as a consequence a proof of the Poincar{\'e}  conjecture.  
 
A few definitions will allow us to more   easily discuss Maskit's work. A {\it Kleinian group} $G$ is a discrete subgroup of $PSL(2,\mathbb{C})$, regarded as the M\"{o}bius transformations acting on the Riemann sphere.  The quotient by $G$  of that portion of the Riemann sphere on which $G$ acts properly discontinuously -- the {\it regular set}  $ \Omega(G)$ -- if nonempty, is a collection of Riemann surfaces.
 
 The {\it limit set} of $G$, $\Lambda (G)$, is the complement of $\Omega(G)$ in $\hat{\mathbb{C}}$.  For any  non-virtually abelian group $G$, $\Lambda(G)$ is the set of accumulation points of the orbit of a point under the action of $G,$ and is a closed perfect subset of $\hat{\mathbb{C}}.$ 
%If the region of discontinuity is nonempty, then the quotient, $\Omega(G)/G,$ of the region of discontinuity by $G$ is a disjoint union of Riemann surfaces.  
The older terminology -- which we shall adopt going forward -- requires that the group $G$ must have a nonempty regular set before it is called Kleinian.  In modern terminology, these are called {\it Kleinian groups of the second kind}. It is largely these groups that Maskit's work is concerned with.  For ease of presentation, we assume going forward that our Kleinian groups are torsion free (that is, have no elliptic elements).

 A {\it Fuchsian} group $G$ is conjugate in $PSL(2,\mathbb{C})$ to a discrete subgroup of $PSL(2,\mathbb{R})$, and the latter acts on the upper half-plane $\mathbb{H}^2$ in the Riemann sphere.  The quotient of the upper half-plane by a Fuchsian group is a single Riemann surface.  There is an extensive theory of Fuchsian groups, which may be regarded as the analytic theory of Riemann surfaces.
 
 Bernie's contributions in the areas of Kleinian groups, low dimensional topology, complex analysis, and geometry of hyperbolic manifolds, include  the planarity theorem (\cite{MR172252}), the study of Schottky groups (\cite{MR220929}),  the Klein-Maskit combination theorems ((\cites{MR192047,  MR223570, MR289768}), and 
 the Poincar{\'e} polyhedron theorem (\cite{MR297997}).   Maskit also was concerned with moduli (deformation) spaces of Riemann surfaces, including Teichm\"{u}ller space, Schottky space, and more generally the deformation theory of Kleinian groups.  
  Maskit's work also included notable results on a  variety of Kleinian groups:  Koebe groups, degenerate groups, component subgroups, Schottky groups and function groups.
By making use of his planarity theorem, Maskit's early work included three counterexamples to several reductions of the Poincar\'{e} conjecture to other conjectures that were proposed by Papakyriakopoulos in the early 1960s (\cite{MR148039}, \cite{MR148040}, \cite{MR145496}).  The planarity theorem has been  of recent use in describing groups of homeomorphisms of  planar covers of the fundamental groups of certain orbifolds and geometrically finite Kleinian groups.  Indeed, the continuing utility of the planarity theorem in this context is evinced in there being a recent new proof by Bowditch \cite{MR4420863}.

 A main thrust of Maskit's career was to understand how one can construct Kleinian (or Fuchsian) groups in different ways that produce collections of Riemann surfaces. 
 In order to carry out this program, Maskit looked not only at the action of a given Kleinian group $G$ on the regular set, but also at the structure of the limit set $\Lambda(G)$, and the action of $G$ on $\Lambda(G).$  Techniques used to generate discrete groups include the combination theorems, the Poincar\'{e} polyhedron theorem, and constructions of Schottky groups.  
Non-Fuchsian examples of Kleinian groups include {\it quasi-Fuchsian groups} (quasiconformal deformations of Fuchsian groups)  and Schottky groups. 
   
Below, we discuss Maskit's  work on Schottky groups, combination theorems, the Poincar\'{e} polyhedron theorem, and briefly, joint work with Alan Beardon.

 If $A_i$ and $B_i$, $i=1,\ldots ,g$, are $2g$ disjoint Jordan curves in the complex sphere that bound a single connected set $D$, and where the exterior of $A_i$ (the component containing $D$) is mapped onto the interior of $B_i$ by a M\"{o}bius transformation $h_i$, $i=1,\ldots ,g$, respectively, then the group $G$ generated by $h_1, \ldots , h_g$ is Kleinian with nonempty regular set, is free on these $g$ generators, and has $D$ as a fundamental domain.  Moreover, all of the non-identity elements of $G$ are loxodromic, and when $g \geq 2$, the 
limit set of $G$ is a Cantor set.  A Kleinian group that may be constructed in this fashion is a {\it Schottky group} and $\Omega (G) /G$ is a genus $g$ Riemann surface (and $(\mathbb{H}^3 \cup \Omega (G))/G$ is a genus $g$ handlebody). Conversely, Maskit showed that a finitely generated subgroup of $PSL(2,\mathbb{C})$
%$M\ddot{o}b_{\mathbb{C}}$  
that is free, purely loxodromic and discontinuous somewhere on $\hat{\mathbb{C}}$, is a Schottky group.  If the generators of the Schottky group for a particular topological surface are taken as parameters, one can form a moduli space -- called the Schottky space --  of Schottky  groups.  The Teichm\"{u}ller space of the underlying surface forms a regular cover of the Schottky space
and  Schottky space is a non-regular cover of moduli space.  Maskit made a number of foundational contributions in understanding the structure of  Schottky space.

The combination theorems are a collection of enhancements of the Klein combination theorem in which conditions are imposed on two Kleinian groups so that the group generated by both of them is also Kleinian, and algebraically is the free product of the groups.  These enhancements include the cases of free products with amalgamation, and HNN extensions.   Maskit's program was to show that all Kleinian groups may be built by using the combination theorems.  This approach is analogous to the approach of Haken and Waldhausen in which a hyperbolic 3-manifold may be cut up into a union of balls and by reversing the steps, the manifold may be reconstructed.  (See \cite{MR1567847}.)   While conditions necessary to apply the combination theorems could not be easily satisfied for the Kleinian groups associated to certain three manifolds,  Maskit's program was almost entirely successful and resulted in a large \lq\lq zoo\rq\rq of cases. (See chapters VIII - X of \cite{MR959135}.) 
For example,  the approach was successful with  {\it function groups}, that is, a Kleinian group having an invariant connected component of its  regular set.  Maskit proved that function groups can be built  from elementary, quasifuchsian, and degenerate groups by using combination theorems a finite number of times, where amalgamations take  place over trivial  or parabolic cyclic groups. For example, a Schottky group is 
a combination of elementary (in this case, infinite cyclic loxodromic) groups.

The upper half-plane $\mathbb{H}^2$ may be endowed with a metric of constant negative curvature ($ds={|dz|\over \mathrm{Im\;}z}$) on which a Fuchsian group acts as isometries.  The Poincar\'{e} polyhedron theorem provides sufficient conditions for a polygon $P$ (with geodesic sides) in $\mathbb{H}^2$, along with isometries of $\mathbb{H}^2$ that pair the sides of $P$, to generate a Fuchsian group for which $P$  is a fundamental region.  Poincar\'{e}'s proof was incomplete, and the Kleinian and Fuchsian groups literature was rife with a number of partial or inaccurate proofs.  Maskit's proof (in \cite{MR297997}), which he averred was the result of conversations he had with Ahlfors, Bers, Magnus, and McMillan, is both correct and more general than Poincar\'{e}'s. Poincar\'{e} and Maskit also provide a generalization of this result for polyhedra in $\mathbb{H}^3$ (see below) in which the group generated by the side pairing transformations is Kleinian.  (See the review of this paper in MathSciNet by J. Lehner.)

 %Briefly, let $G_1$ and $G_2$, be two Kleinian groups, each with nonempty regular sets, that have either only the identity element in common or a cyclic subgroup $H$.  Klein's combination theorem provided  sufficient conditions for concluding that $G_1*G_2$ is also a discontinuous group.  In a series of papers, Maskit provided a number of enhancements to the basic theorem that allowed constructions, and descriptions, of Kleinian groups of a wide variety of sorts.  (FIX THIS) 
 
% Given a surface, $S$, the planarity theorem provided a representation of each regular covering surfaces of $S$ that are planar by a collection of loops on the surface $S$. 
 
 %%%%%%%%%FIGURE%%%%%%%%%%%%%%%%
\begin{figure}[htb]
 \begin{center}
     \includegraphics[width=3.0in]{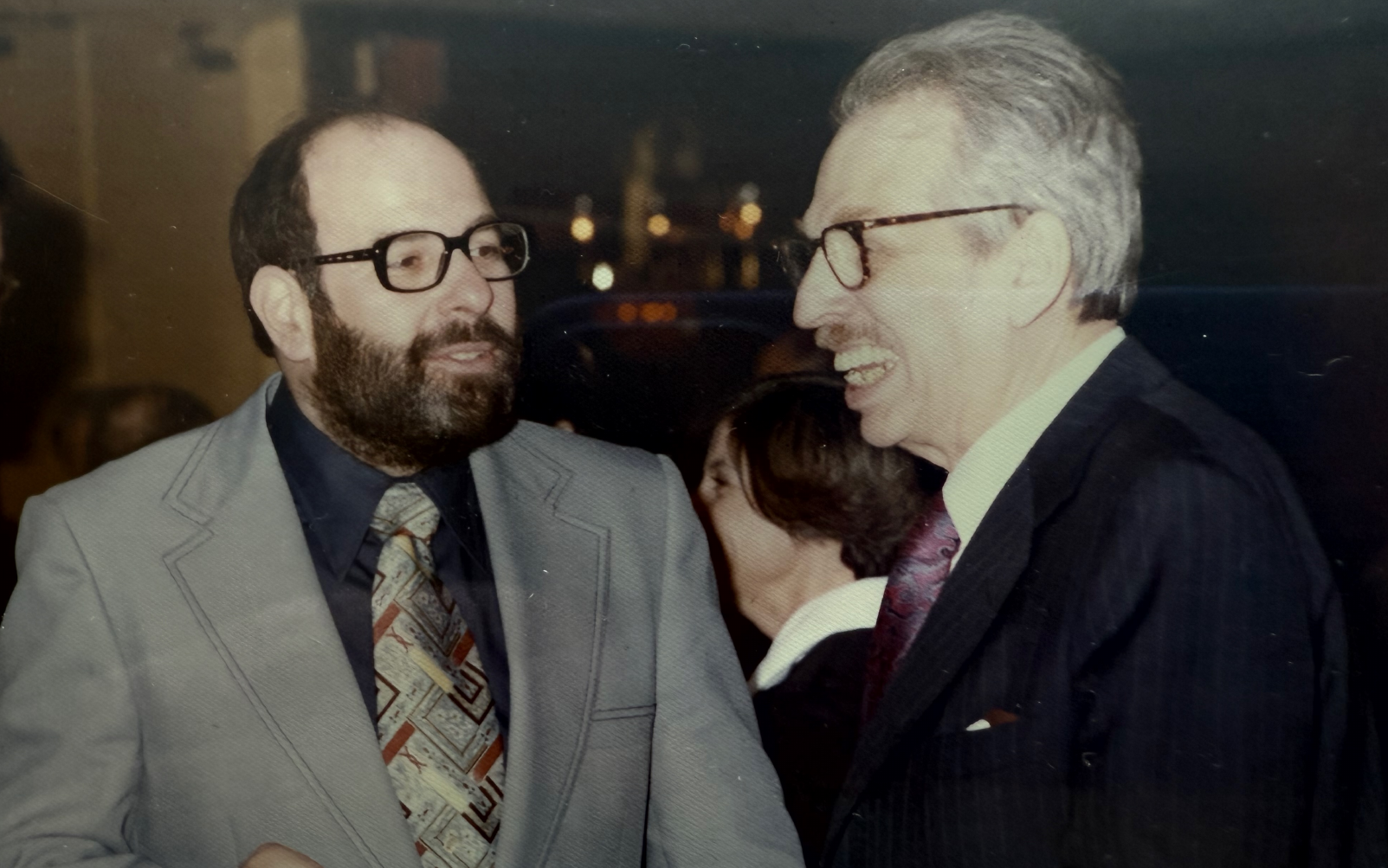}
     \vspace{-10pt}
    \caption{Bernie and  Lipman Bers}
    \vspace*{-.2in}
    \label{}
  \end{center}
\end{figure}
%%%%%%%%%%%Figure%%%%%%%%%%%%

     In 1883, in an attempt to treat the theory of Kleinian groups in roughly the same manner as the theory of Fuchsian groups, Poincar\'{e} extended the action of the  M\"{o}bius transformations from the Riemann  sphere to the upper half-space $\mathbb{H}^3 = \{(z,t)\ |\ z\in \mathbb{C}, t>0\}$,  with boundary 
     $\partial\mathbb{H}^3=\mathbb{C}\cup \{\infty\}=
     \widehat{\mathbb{C}}.$  The M\"{o}bius transformations now act as isometries on $\mathbb{H}^3$ with metric 
     $ds^{2}={|dz|^2+dt^2 \over t^{2}}$, and as with Fuchsian groups, it is elementary to show, as a consequence, that a Kleinian (discrete) group must act properly discontinuously on $\mathbb{H}^3$.   For a Kleinian group $G$ acting on $\mathbb{H}^3  \cup \Omega(G)$, where $\Omega(G) \ne \emptyset$ (and recalling that we are assuming for the rest of this discussion that   $G$ has no torsion) the quotient, $\Omega(G)/G$, which was earlier regarded as a disjoint union of Riemann surfaces, may now be regarded as the collection of boundary components of the hyperbolic 3-manifold (as opposed to orbifold, if there is torsion) with boundary, $(\mathbb{H}^3  \cup \Omega(G))/G.$

    After Poincar\'{e}'s effort, the theory of Kleinian groups languished somewhat for a number of decades, but efforts to develop the theory picked up in 1964 with Ahlfors' finiteness theorem  (\cite{MR167618}). Maskit's contributions to Kleinian groups, beginning at that time, were generally of a geometric/topological nature and, as already discussed, included the  Poincar\'{e}'s polyhedron theorem, the Klein-Maskit combination theorems and the planarity theorem.   
   In particular, Maskit developed the theory of certain special Kleinian groups that he was able to build by making use of the Klein-Maskit combination theorems.  These were the B-groups and degenerate groups.  In the late 1960s, Greenberg (\cite{MR200446}) was able to use the existence of {\it degenerate groups} (a function  group whose limit set contains more than two points and  whose regular set is connected and  simply connected) in order to prove that there are finitely generated Kleinian groups in which, for each,  there is not associated a finite-sided fundamental polyhedron (see below) in $\mathbb{H}^3$ for the action of the group on $\mathbb{H}^3$.

 Maskit's work in part motivated the work of Albert Marden in (\cite{MR349992}).   To quote Marden,  
%Nearly concurrently, and  partly motivated by the above,  Albert Marden's aspiration in  \cite{MR349992} was to realize Poincar\'e's original hope of treating the theory of Kleinian groups acting on hyperbolic space  in much the same way as the theory of Fuchsian groups.    To quote Marden,   
    
    \begin{quote}
     In contrast to the analytic approach of Ahlfors and Bers, Maskit has pioneered in the study of Kleinian groups by purely geometric methods in the complex plane. His work complements the insight provided by the analytic approach and yields some very deep and fundamental knowledge concerning certain special classes of groups. In particular Maskit develops the concept that the class of \lq nice\rq  \ groups, that is the class for which generalizations of the classical results for Fuchsian groups can be fruitfully sought, is the class of constructible groups. These are the groups that arise from his far-reaching generalizations of the Klein combination theorems. Although we have taken a different path in this paper, Maskit's work has been a great influence. Indeed, there are many points of contact between his work and ours although the precise relation remains unclear.  
     \end{quote}

 Marden's program was to realize Poincar\'e's original hope of treating the theory of Kleinian groups acting on $\mathbb{H}^3$ in much the same way as the theory of Fuchsian groups, the latter being a way of considering the theory of Riemann surfaces in which a Riemann surface $S = \mathbb{H}^2/G$ is considered as a quotient of $\mathbb{H}^2$ by a Fuchsian group $G.$  Indeed later, Thurston made much use of Marden's work (and the work of Haken and Waldhausen) in developing the theory of hyperbolic 3-manifolds, paving the way to Thurston's geometrization conjecture, the proof of which was completed by Perelman in 2003.  Analogous to the classification of Riemann surfaces a century earlier, which states that every Riemann surface admits a metric of constant curvature, the geometrization theorem states that any oriented irreducible closed   3-manifold  may be decomposed into pieces, each of which has one of eight unique geometric structures of finite volume.  (Here, irreducible means that every embedded 2-sphere bounds an embedded 3-ball.) Among those eight, the hyperbolic structure is the principal one, and most of the work discussed herein, including Maskit's, is concerned with this particular geometry.

Although developments that led to the geometrization conjecture are often described in terms of surgical techniques on 3-manifolds, it is fair to say that Maskit's combination theorems provide part of the technology that was used in formulating this most fundamental result of the late 20th and early 21st century.  Among other notable consequences, the geometrization theorem implies the Poincar\'{e} conjecture.  It is also striking that Maskit's early work on the planarity  theorem precluded certain approaches to the Poincar\'{e} conjecture -- see Keen's contribution below -- but his combination theorems contributed to developments that eventually led to a proof.

 Many of the types of Kleinian groups considered by Bernie were not just finitely generated, but also satisfied the stronger condition of being {\it geometrically finite}.  A Kleinian group $G$ is {\it geometrically finite} if there is a finite-sided fundamental polyhedron for the action of $G$ on $\mathbb{H}^3.$  For such groups the 3-manifold with boundary, $M=(\mathbb{H}^3 \cup \Omega(G))/G$, may be understood as being \lq\lq nearly" compact, that is, $M$ is compact except for a finite collection of \lq\lq cusps" that arise in one-to-one correspondence with conjugacy classes of Maximal parabolic subgroups in $G$.   The behavior of geodesics in the interior of $M$ may also be understood by examining the properties of the limit set of the geometrically finite group $G$.  Indeed, the work of  Ahlfors, Beardon, Greenberg, Maskit, Marden, Thurston and others resulted in a number of equivalent analytic, geometric, topological, and dynamical \lq\lq finiteness conditions" that hold for a given geometrically finite group $G.$ These characterizations of geometric finiteness include the well-known $thick-thin$ decomposition of the quotient manifold $ M = \mathbb{H}^3/G,$  given by Thurston.

   Though Maskit's work generally focused on the 2-dimensional theory, in Maskit's most cited paper, joint with Alan Beardon  (\cite{MR333164}), the authors provide a characterization of a geometrically finite Kleinian group $G$ entirely in terms of the characteristics of the limit set of $G$.  Indeed, $G$ is geometrically finite if and only if the limit set of $G$ consists of {\it points of approximation} (called {\it conical limit points} in more recent literature), and
 {\it (cusped) parabolic fixed points} (more recently called {\it bounded parabolic fixed points}). This characterization of the limit set of a geometrically finite group $G$ is part of a collection of five closely related equivalent notions of geometric finiteness, and all have been generalized in a number of ways. For discrete groups of isometries of $\mathbb{H}^n$, Bowditch, Tukia and others generalized the notion of geometric finiteness (\cite{MR1218098}), and for discrete groups of isometries of $Isom(X)$ where $X$ is an $n$-dimensional {\it pinched Hadamard manifold}, that is, $X$ is a complete, simply connected Riemannian manifold with all
sectional curvatures bounded between two negative constants, Bowditch has carried out the generalization (\cite{MR1317633}).  Other generalizations include convergence groups.  Conical limit points are related to the behavior of geodesics in the quotient manifold for discrete groups of isometries in all of these settings, where a limit point $z$ is a conical limit point if a geodesic ray ending at $z$ projects to a geodesic in the quotient manifold that returns infinitely often to a compact set.  Consequently, returning to the category of Fuchsian groups, there are rather interesting connections between the behavior of geodesics on the quotient of $\mathbb{H}^2$  by the modular group $PSL(2, \mathbb{Z}))$ to number theory. 
 For geometrically finite manifolds of infinite volume, the proof of the ending lamination conjecture follows from the work of Maskit, Kra and others.

 Any expert in Kleinian groups will be familiar with the work of Bernard Maskit, but possibly anyone who has even considered study in the area will be familiar with Maskit's 1988 book,  {\it Kleinian groups} (\cite{MR959135}).  In it he covers many essential topics -- indeed, many of the topics chosen are those that Maskit had a large hand in developing.  Although pains are taken in many cases to approach topics in a way that allows for some generality (for example, discrete subgroups of the isometries of hyperbolic $n-$space are considered, rather than isometries of hyperbolic $2-$ or $3-$space), the approach is largely 2-dimensional. That is, the approach is largely oriented toward developing the planar theory of Kleinian groups by exploring their action on the Riemann  sphere, rather than the 3-dimensional theory of hyperbolic manifolds, which results by looking at the action of a Kleinian group on the upper half-space.  One should observe that this approach may make it appear that Maskit's work is somewhat removed from the exciting work on geometrization and hyperbolic 3-manifolds that has occurred over the last 40 years, but this is mostly an issue of point of view, since these areas are either inextricably intertwined or in some cases identical.
 
 %%%%%%%%%FIGURE%%%%%%%%%%%%%%%%
\begin{figure}[htb]
 \begin{center}
     \includegraphics[width=3.5in]{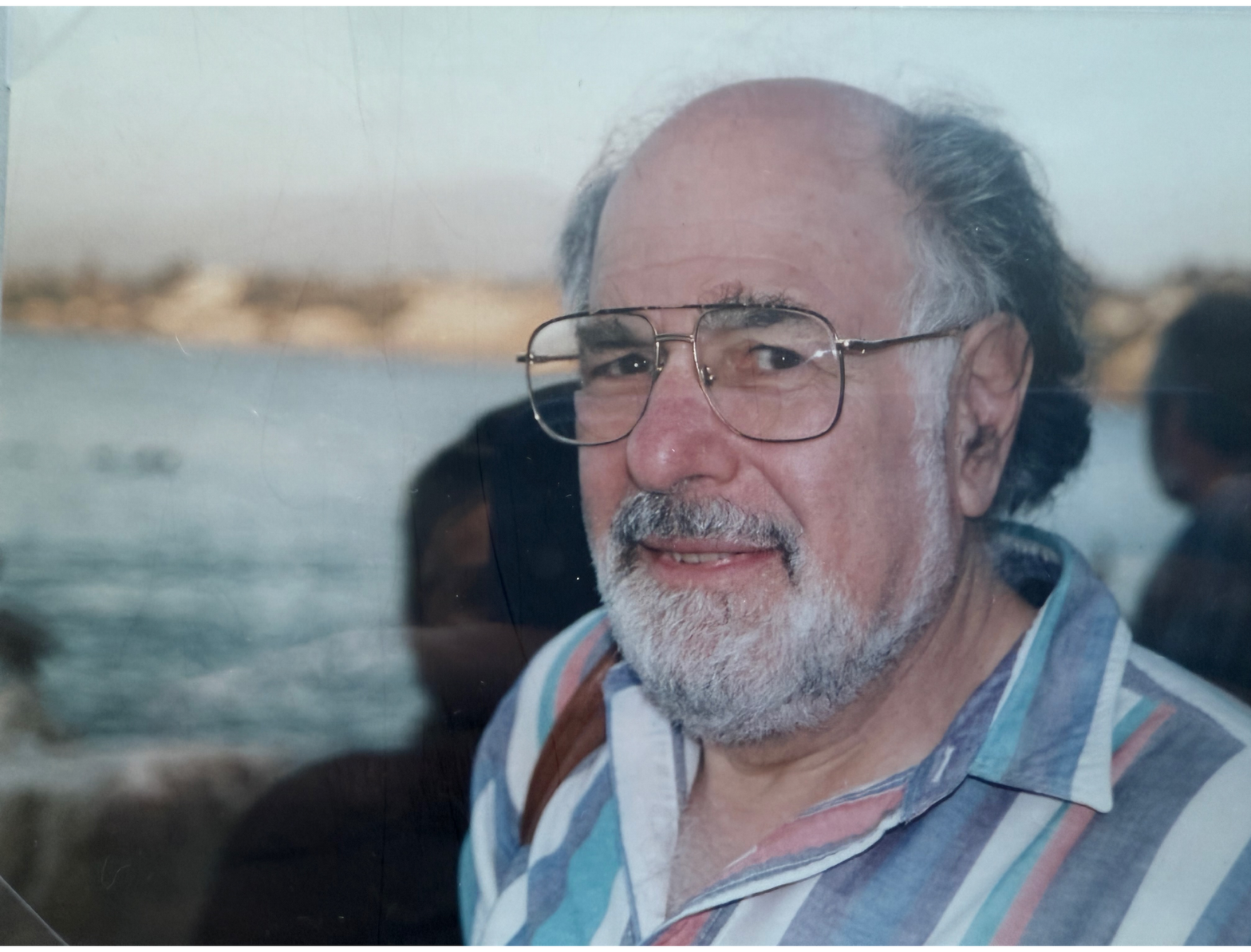}
     \vspace{-20pt}
    \caption{Bernie  in Chile}
    \vspace*{-.2in}
    \label{}
  \end{center}
\end{figure}
\begin{figure}[htb]
 \begin{center}
     \includegraphics[width=3.0in]{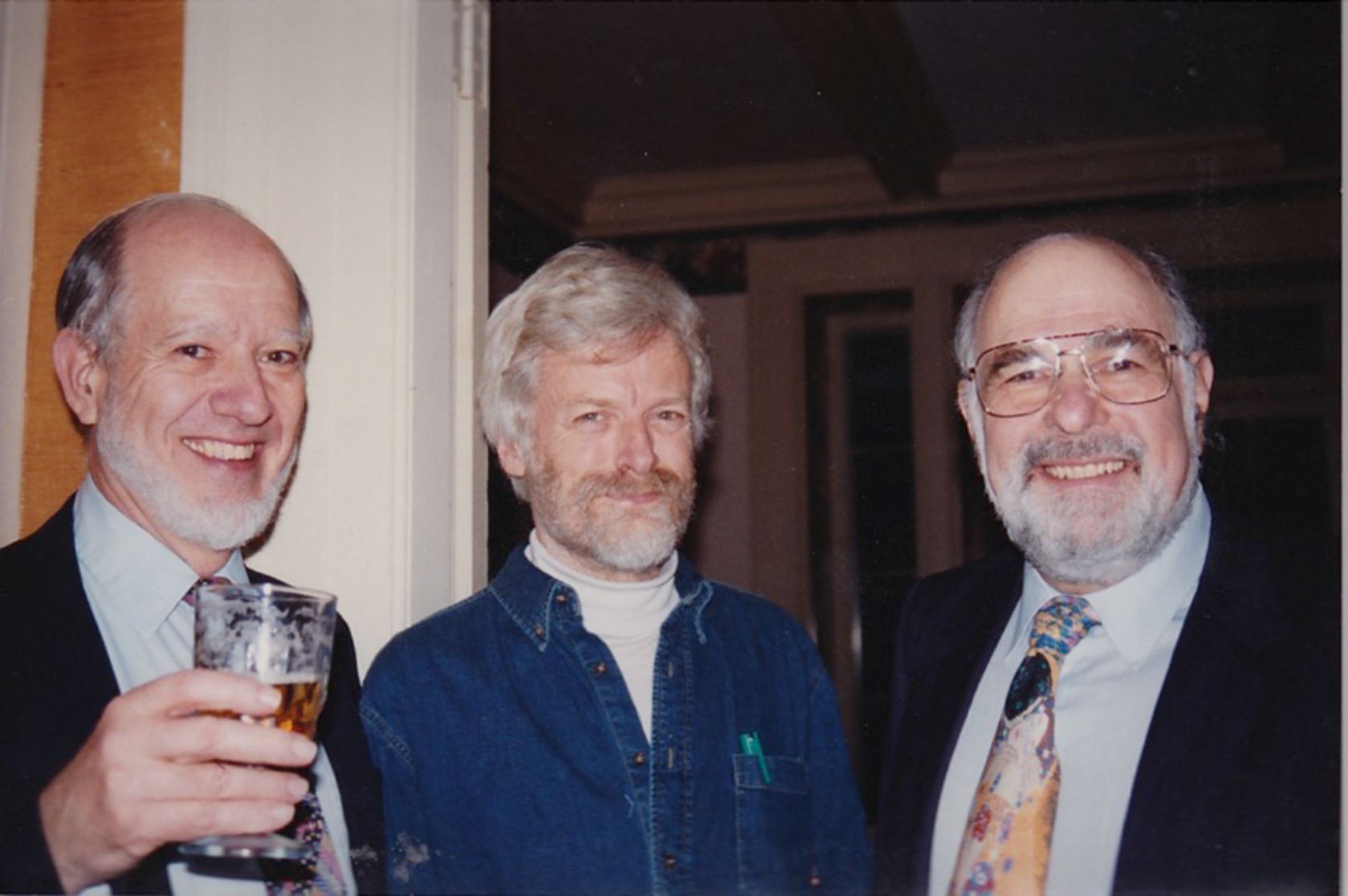}
     \vspace{-10pt}
    \caption{Clifford Earle, William Harvey and  Bernie}
    \vspace*{-.3in}
    \label{}
  \end{center}
\end{figure}
%%%%%%%%%%%Figure%%%%%%%%%%%%

      %%%%%%%%%FIGURE%%%%%%%%%%%%%%%%
\begin{figure}[htb]
 \begin{center}
     \includegraphics[width=3.0in]{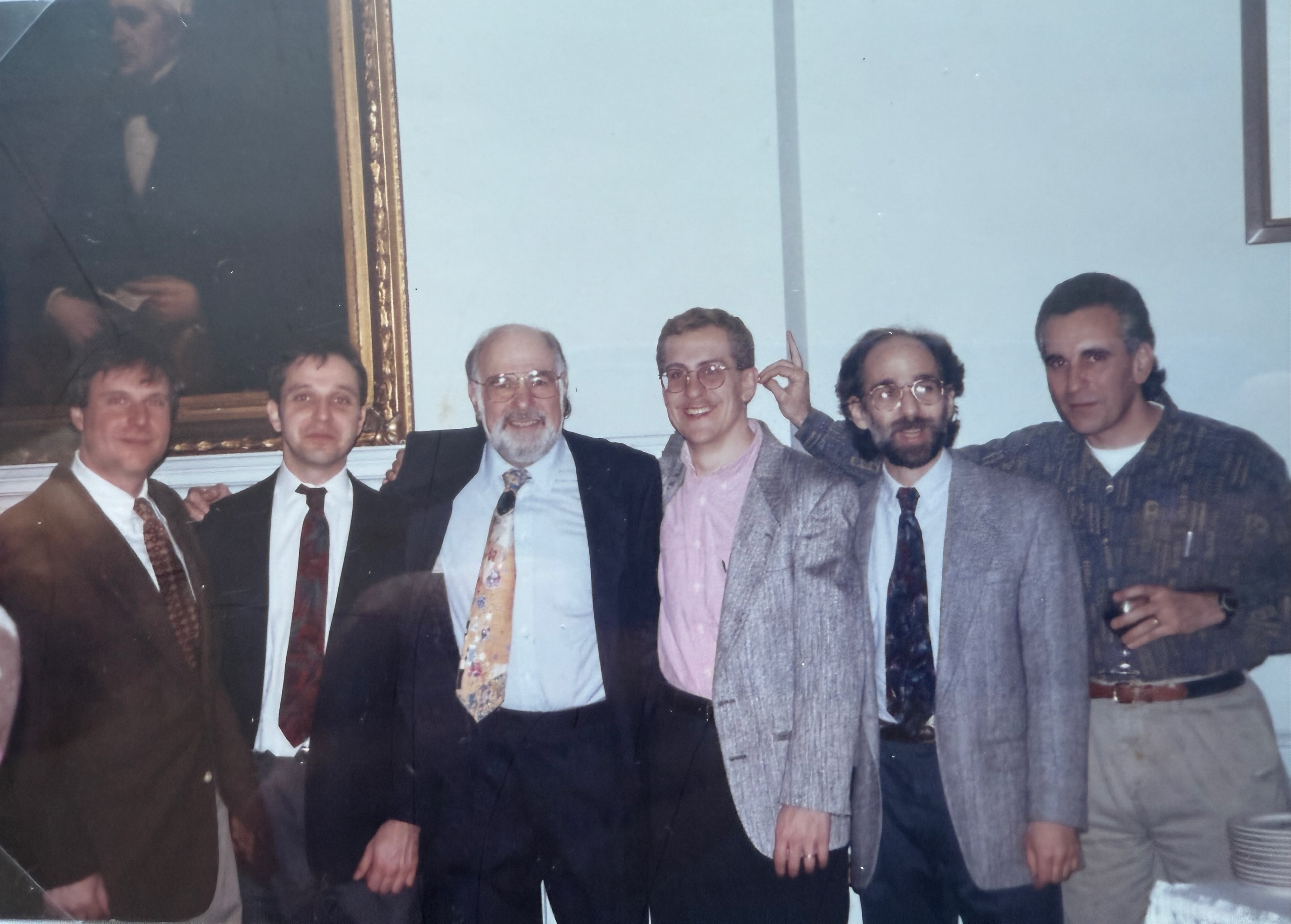}
     \vspace{0pt}
    \caption{Bernie's 60th birthday with students (from left to right) Perry Susskind, Ed Taylor, Bernie, Jim Anderson, Andy Haas and Ara Basmajian}
    \vspace*{-0.2in}
    \label{}
  \end{center}
\end{figure}
%%%%%%%%%%%Figure%%%%%%%%%%%%

  %%%%%%%%%FIGURE%%%%%%%%%%%%%%%%
\begin{figure}[htb]
 \begin{center}
     \includegraphics[width=3.0in]{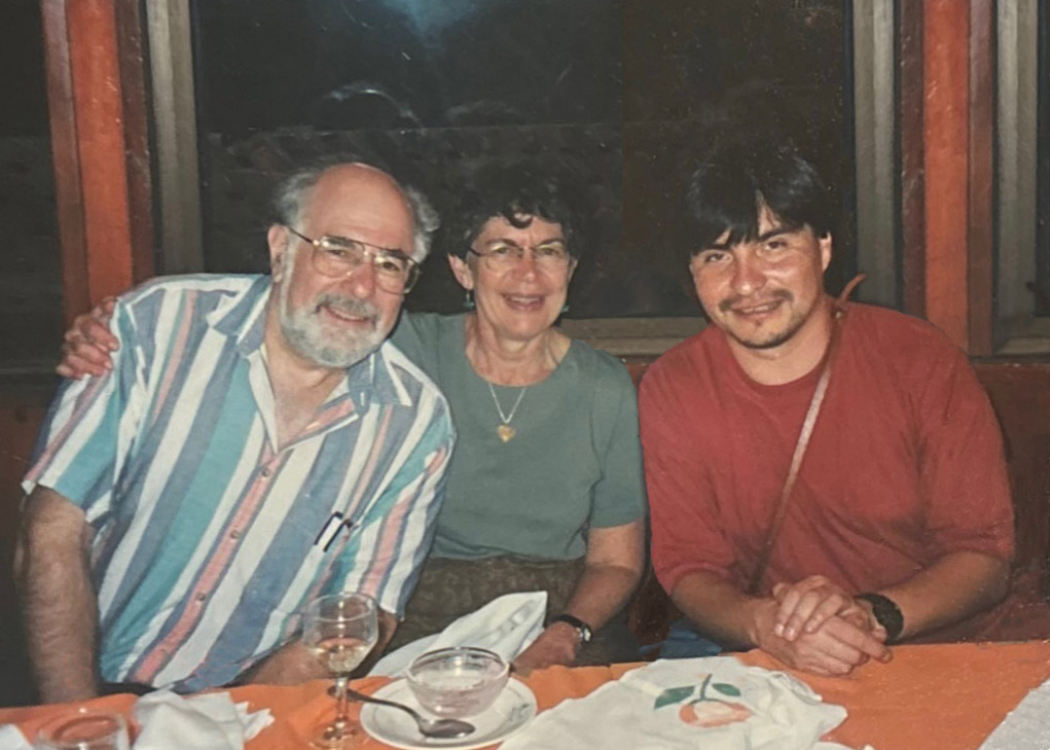}
     \vspace{-10pt}
    \caption{Bernie, Wilma, and Ruben Hidalgo}
    \vspace*{-.3in}
    \label{}
  \end{center}
\end{figure}
%%%%%%%%%%%Figure%%%%%%%%%%%%

Each of his students had their own experience of Bernie as an advisor, though all would agree  that Bernie was a most important figure for his students.  While Bernie did not provide an over-abundance of advice and guidance to his students, he  provided just what was necessary for each, and consequently allowed each of his students, according to their individual abilities and mathematical maturity, to stand on their own two feet.  The largely topological arguments of many of Bernie's papers required highly demanding precision,  rigor and meticulous attention.  Beyond just the shepherding of students through the usual process of turning a particular but possibly vague mathematical problem or idea into a viable thesis, Bernie quietly and modestly modeled by example the kind of meticulous and careful thinking that guided all of us in our subsequent mathematical (or not mathematical) careers.  

When stating theorems in his Kleinian groups course, Bernie would often write just the conclusion, and announce that he would be filling in the hypotheses as needed during the proof, somehow mimicking a creative process more engaging than just \lq\lq theorem--proof."   Peter Buser related a rather interesting criterion that Bernie might employ in order to deem a paper finished. \lq\lq It was in Oberwolfach, or perhaps at some other place. I was sitting at a table with Bernie and we were discussing publishing mathematics. My concern was a paper that was done, but yet not finished, and I couldn't get it done. \lq Hmm', he said, and after a while: \lq and then you improve a passage and change a notation, and later you feel that the earlier version was better, and you think it over again. You know, when that happens, it's the sign that the material has matured.' I often remember this scene when I'm working on a text."  

One of us remarked at a conference banquet held in his honor at NYU in March of 2008  that when teaching classes, Bernie's students benefited from, and appreciated watching, listening, and learning, from  instances in which Bernie would think through a nontrivial argument on the fly without having prepared it in advance.   (As a response, Bernie quipped that he did not realize he was being a good teacher; rather he claimed he thought he was just being lazy about preparing his classes in advance.)  At times, after making an assertion in class, Bernie would ask, \lq\lq Do I have to worry?" Sometimes the answer was \lq\lq Yes," and significant explication would follow.  Other times, there ensued a real-life instance of the apocryphal story of the student challenging the professor's statement that some observation made in class  is trivial, wherein the professor stands in silence for several minutes engaged in deep thought and then reiterates, \lq\lq Yes, it's trivial." In this instance, Bernie would ask, \lq\lq Do I have to worry,"  stand for a minute or two in silence,  thinking deeply, and then proclaim, \lq\lq No, I don't have to worry."  His students were usually relieved that he -- and we -- didn't have to worry.  

On a personal level, Bernie was reserved with his students but would invite us to gatherings at his house where there was often excellent food and drink.  Lastly, we mention a final chapter in Bernie's intellectual life.  A few years before he retired, Bernie became interested in an area of computational psycholinguistics called {\it multiple code theory}, developed by his wife and intellectual partner Wilma Bucci.  Bernie made significant contributions in this area (see \cite{Su1}) largely by introducing the use of the technique of mathematical smoothing for certain mathematical measures on spoken narratives.    All of us agree, likely due to Wilma's influence, that Bernie became increasingly content as time passed.  They enjoyed good food and drink, and gathering with friends.  Bernie enjoyed tennis for a time -- he was a handball player in his early years -- and few might guess that Bernie and Wilma learned how to dance the Tango.  

Our motivation for writing this memorial article began with the impetus to honor our advisor and make others aware of his work and role in important developments in 2- and 3-dimensional topology and geometry.  Years ago, as new Ph.D.s, we were perhaps only partially equipped to fully understand the breadth and depth of Bernie's contributions in these areas. With the perspective of much time working in these areas ourselves, the process of revisiting Bernie's contributions to mathematics has been a highly rewarding endeavor.

 %%%%%%%%%FIGURE%%%%%%%%%%%%%%%%
\begin{figure}[htb]
 \begin{center}
     \includegraphics[width=3.0in]{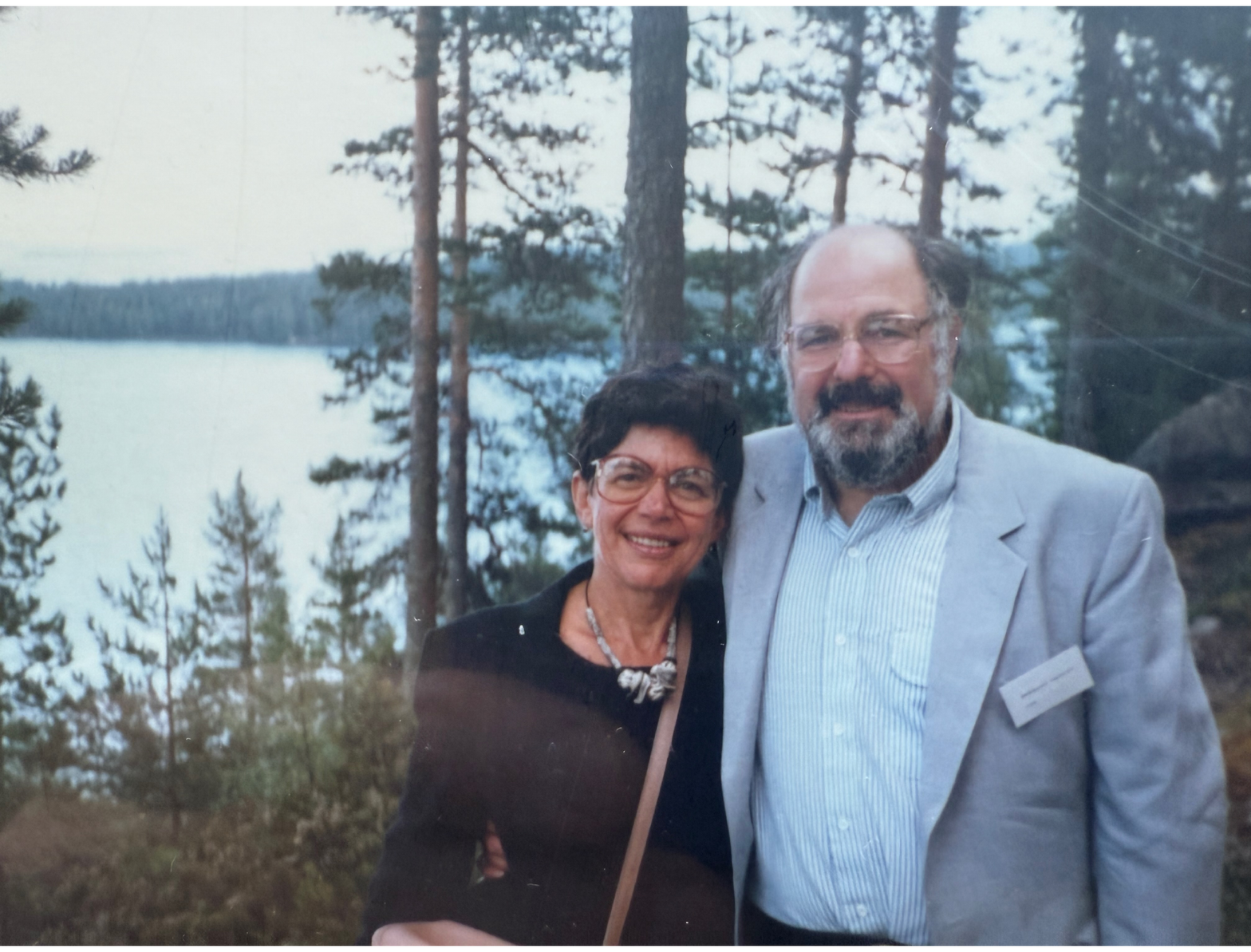}
     \vspace{-10pt}
    \caption{Bernie and Wilma in Chile}
    \vspace*{-0.4in}
    \label{}
  \end{center}
\end{figure}
%%%%%%%%%%%Figure%%%%%%%%%%%%

%\subsection*{William Abikoff}
  \vskip5pt
\noindent\large{\bf William Abikoff}
\vskip5pt
\normalsize
I first met Bernie at the Stony Brook Riemann Surfaces conference in 1969. I was an engineer at Bell Labs and had recently passed mathematics doctoral exams at Brooklyn Poly. I had been talked out, actually laughed out, of my interest in a particular problem about sums of random functions and had come uninvited to that conference. The only attendee I knew was Lesley Sibner. The excitement in the air was overwhelming. After dark, participants talked and socialized while I, in desperation, drove to the Courant Institute library, the only research library open at night, and tried to figure out what all those excited people were talking about.  A few of them not only welcomed me but they made me feel that I belonged. Even now, over a half-century later, I treasure the moments when I recall them --- Bernie, Lipa Bers, Ernie Rauch and Bob Gunning.  It didn't take long for me to become an auxiliary member of the \lq\lq Bers Mafia."

Probably at Lipa's suggestion, I visited Bernie at MIT around Christmas of 1971. We soon bonded both professionally and personally. At that time, the Kleinian groups crowd consisted of his contemporaries. And me. I became something of a kid brother.

We shared non-standard paths from The Bronx into and through mathematics. Bernie
had started undergraduate study at RPI,  then withdrew because of family problems. He then spent some time assisting in his family's business, while completing his undergraduate study at NYU. Starting in his graduate years, his main focus lay in the interplay and consequences of algebraic and geometric constructions. Later on, his curiosity moved him to consider questions raised in the work of his, in so many ways, life partner Wilma Bucci. 

%We shared non-standard paths into and through mathematics. Bernie had started undergraduate study at RPI only to withdraw in order to assist in his family's business. Not long thereafter, he completed undergraduate study.  Starting in his graduate years, his main focus lay in the interplay and consequences of algebraic and geometric constructions. %Much later on, when nearing retirement, his curiosity moved him to consider questions raised in the work of his, in so many ways, life partner Wilma Bucci. 

One time we met to discuss a question we were considering; to be honest, it was a shared annoyance. He had a proof with a gap; I had a proof, the same proof in a different setting, also having a gap. But the gaps were different. Together, we had the proof. 
 
%The central theme of Bernie's work lay in the construction of discrete groups of M\"obius transformations. For the most part, Bernie stated his results in terms of M\"obius transformations acting on the Riemann sphere but, as first observed by Poincar\'e, the M\"obius action on the Riemann sphere extends to the ball it bounds and that action is isometric in the hyperbolic metric on  the ball. 

While much of his work predated the widespread use of computers in geometry, his constructions are  computable. In some cases, Thurston's Bounded Image Theorem is paired with the constructions --- this interdependence is key to most of the proof of Thurston's Hyperbolization Theorem. Only manifolds fibering over the circle require a different approach. While Thurston was developing the ideas for hyperbolization, he lectured quite publicly before he had complete proofs. Thurston used at least one combination theorem that Bernie hadn't yet proved. In one lecture Bill gave at Columbia, Bernie informed him that the combination theorem he used in the lecture was not yet proved.  Felix Klein is credited with the first construction (See the reference in \cite{MR192047}.) In the eighty year period between Klein's work and Bernie's, there were few constructivist contributions to the literature other than in the Fuchsian, i.e. invariant disk, case. The two that come quickly to mind are Lester Ford's 1926 book and the mysterious Fenchel-Nielsen manuscript; the latter was announced as appearing imminently in the late 1940's and was finally published in 2003. 

Bernie's work has much contemporary interest in a variety of fields, including complex analysis, hyperbolic geometry and knot theory. It had reenergized combinatorial group theory and is a spiritual precursor of  geometric group theory. For me to write of  its future influence would require prophetic skills that I do not possess.

%\subsection*{Wilma Bucci}
  \vskip5pt
\noindent\large{\bf Wilma Bucci}
\vskip5pt
\normalsize
I met Bernie on Christmas Eve 1978, almost half a century ago, when he was on a skiing trip with his two younger sons, and staying with his sister who was a psychologist and friend of mine.  Back then, Bernie's three sons all were under the age of 20; my son and daughter were in the same age range. We shared a strong wish to do the best we could for our family.  Bernie was truly passionate about mathematics, but also about his family.  There were many instances of how conflict between family and work played out in our life together.  As an example, about a year after we met, Bernie was involved in organizing a conference focused on Bill Thurston's work, which was new at the time.  Bernie and I, while traveling in Maine, found the perfect venue for the conference at Bowdoin College.  We drove to the meeting  together, which took place in the summer of 1980.  A young relative was left with my children in New Jersey; Bernie's sons were with their mother.  Each day at the meeting, Thurston presented his work.  Each day after his talks, Bernie, Lipman Bers and other colleagues met impromptu in sessions devoted to unraveling what Thurston had said. Bernie and I would then go back to our campsite near the meeting, and before settling down for the evening called home to check on our children. (There were no cell phones in those days so we had to find a pay phone near the site.)  About the third night, I got one of those responses that parents dread: \lq\lq Everything is alright now.\rq\rq    My daughter had been in a bicycle accident, had been taken to the emergency room, and had broken her arm, though she was said to be okay.  I told Bernie I had to go home right away; he said he was going with me. He told Lipa (Lipman Bers) and their colleagues. Lipa said they needed Bernie to stay. Bernie said he was going with me.  Mary Bers said of course he has to go with Wilma.  Bernie drove me home and we dealt with that family problem together, as we did with many problems on many levels throughout the years. 

As the children became more independent, we did extensive  traveling, went together to many conferences, and to visit colleagues, in Europe and South America, including a wonderful year in France, based at IHES, just outside Paris, where we went to market, cooked, and drank wine. We always cooked together, wherever we were, Bernie mostly as prep chef.  Bernie was an experienced hiker and camper.  We drove across the U.S. twice, camping in national parks and forests and staying for the summer in Berkeley where Bernie worked at MSRI and I commuted to UCSF. Bernie played the recorder seriously and well, and we did all kinds of dancing -- folk dancing, contra dancing, all types of ballroom dancing and then became addicted to Argentinian tango.  

 While Bernie and I joined our families and our lives, we went our separate ways in our research for the first several decades we were together.  Around 2000, as he neared retirement from Stony Brook, Bernie moved from a full focus on mathematics to an increasing interest in research questions associated with my work as a psychologist and psycholinguist.  Bernie brought his creative ways of thought to difficult questions associated with my research. When he looked at my work, which involved studying linguistic processes using computer assisted measures, Bernie tended to first characterize the problems as impossible to solve, but he then proceeded to solve them, in the process introducing mathematical techniques that those in my field had not considered.  Bernie continued developing new ideas to the last days of his life and left us, as part of his legacy, a number of ideas for projects that he planned to pursue.  Some of these new ideas are presented in the last paper Bernie wrote, to be published in a special issue of the Journal of Psycholinguistic Research honoring his work.  Bernie also guided students working in this area on dissertations, papers and posters.  Research centers in the U.S., Canada, Italy, Israel, and Argentina apply Bernie's measures of the referential process, and work by colleagues in these centers continues to extend the theory and applications.
% \subsection*{Fred Gardiner: Some memories} 
    \vskip5pt
\noindent\large{\bf Fred Gardiner: Some memories}
\vskip5pt
\normalsize
  Bernie leaves an extensive legacy of publications and a  book that has become an essential reference in Kleinian groups.  But here I also want to focus on his  generosity, fairness and ability to come up with pointed comments.  These qualities enabled him to draw together three generations of successful working mathematicians. It started in the early 1960's when he earned his Ph.D. and began a continuing collaboration with his thesis adviser, Lipman Bers, at NYU over many years. 
  %Over his long career, Bernie worked in the areas of analysis, geometry and topology, particularly on Riemann surfaces, Kleinian groups, and their moduli spaces.   Looking back now, his success seems  to have been guaranteed from the beginning by the choice of this  beautiful subject.   
Riemann surface theory had become of interest in the 1930s to nuclear physicists.  In particular 
 Hermann Weyl had ten years before published his monograph ``Die Idee der Riemannschen Fl\"{a}che."  This was coupled with Schr\"odinger's earlier use of complex analysis to give discrete solutions to the wave equation which described the discrete possibility for radii of orbitals in an atom.
%Later on in the 1950's and 60's all sorts of related topics emerged including computer experimentation, fast computation and ability to save and retrieve massive amounts data. 
%Today, the internet itself seems (I think falsely) to have stolen first prize for the greatest invention of all time.
%In early times the wheel had  the blue ribbon and a bit later some thoughtful souls suggested the alphabet.  More recently,  a popular  comedien, pretending  to be thousand years old, said it was saran wrap and Cole Porter in one of his songs suggests cellophane. 
But, in  Bernie's life it was the richness of the subject  of  Kleinian groups and Riemann surfaces that drew so much attention.  
 In my memory several special moments stand out. 
 The first occurred at  Bers' complex analysis seminar, in which Bernie was an active participant,    on the fourth floor of the Columbia University Math Building. Dennis Sullivan, having just returned  from IMPA in S\~ao Paulo, came to speak about a paper he wrote jointly with Ma\~{n}\'{e} and Sad.  That paper was about iteration of rational functions, and in the middle of working on it their paper contained the  idea of a holomorphic motion.
  It was a way to conclude that a certain curve of deformations of a dynamical system acting injectively on a limit set would necessarily be realized by a quasiconformal conjugacy.   
That turned out to be the beginning of the studies of holomorphic motions. 
%Slodkowski's extension theorem and the importance of looking at holomorphic functions of two complex variables, which remains a fundamental idea.  
But at roughly the same time, there was what became known as {\it Sullivan's dictionary}, a collection of correspondences that tied the ideas and theorems of Kleinian groups, including the work of Ahlfors, Bers, Maskit and others, with the burgeoning subject of dynamical systems and the iterations of rational maps.  

Another moment happened in the mid-1980s when Bill Thurston came to the same seminar  to explain parts of his extraordinary theory of measured foliations and  laminations, which had begun in the mid 1970s publication of  ``Travaux de Thurston sur les surfaces."  All the while Bernie's expertise and pointed questions contributed immensely to this environment. 
      
When Bernie was at a seminar you would  find other mathematicians working seriously on related topics about which he was an expert and you would find collaborators.   The wonderful Chapter VIII of his book on Kleinian groups is titled ``A Trip to the  Zoo."  That's where he took us and where we had a great time.
%\subsection*{Jane Gilman} 
  \vskip10pt
\noindent\large{\bf Jane Gilman}
\vskip5pt
\normalsize
I met Bernie when I was a graduate student at Columbia and he was already a well established mathematician. We were both students of Lipman Bers, but Bernie was reputed to be the best of Bers' many students working on Teichm\"uller theory. I was in awe of Bernie, but later, when I collaborated with him, I learned that one did not have to be in awe of Bernie, just of his mathematical powers. He was very approachable and very nice. I also got to know some of Bernie's students, especially those from his 1982-95 Kleinian groups--Stony Brook period. In addition to publishing important papers in hyperbolic geometry and Kleinian groups, these students were all very nice and have been a welcome addition to my mathematical world.
My collaboration with Bernie began at an MSRI semester on Kleinian groups and involved the two generator $PSL(2,\mathbb{R})$ discreteness problem: given a pair of hyperbolic elements with disjoint axes, determine whether or not they generate a discrete non-elementary group. This problem has a long history of publications with incorrect or incomplete solutions. Perhaps because the problem is easy to state, people expected a simple closed form solution. Our joint 1991 paper gave a complete geometric algorithm for determining discreteness or non-discreteness. This one collaboration influenced much of my subsequent work on two generator groups.
I remember Bernie and his wife Wilma on an overnight trip to Leningrad from a conference at the University of Joensuu in Finland in 1987. At the time laptop computers were rare, but Bernie had one and brought it with him on the trip. Soviet customs officers did not know what to make of it. We were delayed for hours. Eventually the laptop was allowed to enter the Soviet Union.  
%\subsection*{Linda Keen: Remembrances of Bernie Maskit}
  \vskip5pt
\noindent\large{\bf Linda Keen: Remembrances of Bernie Maskit}
\vskip5pt
\normalsize
  Bernie and I were students of Lipman Bers in the early 1960s at NYU. It was an exciting time. Ahlfors and Bers had just proved a major theorem that made it possible to put a holomorphic structure on the moduli spaces of Riemann surfaces and Kleinian groups, and the structure of Kleinian groups was something into which Bernie had tremendous insight and understanding.
Beginning with his thesis, he became the world's expert on the intricacies of Kleinian groups, for example, how complicated the limit sets could be. This was particularly important in studying what happens as one moves to the boundary in these spaces.
After finishing his thesis, Bernie spent two years at the Institute for Advanced Study in Princeton, and I joined him there the second year. That was a time when it looked very promising that one could prove the Poincar\'e conjecture for manifolds of dimension less than 5. It was the holy grail then, and attracted the best and the brightest, so it was a not little thing that a result of Bernie's, his Planarity theorem, showed that one of the main approaches that would have proved the conjecture in dimension 3 was false. It took another thirty-eight years before the conjecture was proved.
Cliff Earle was also at the Institute that year and the three of us had a wonderful time talking about math together. We read one another's work. I'm particularly grateful to Bernie for carefully reading my paper on moduli for Fuchsian groups and helping clarify its exposition. This was also the year of the first conference on Kleinian groups organized by Ahlfors and Bers in which they showcased their students. It was in New Orleans where we not only talked math, but sampled the food and wine ordered by Ahlfors at many of the famous restaurants. It was the first of many such conferences which grew as the numbers of their students multiplied. Bernie was always front and center in this community that launched so many careers.
For over sixty years, Bernie was a good friend and excellent colleague.
%\subsection*{Irwin Kra: The day after and memories}   
\vskip5pt
\noindent\large{\bf Irwin Kra: The day after and memories}
\vskip5pt
\normalsize
Yesterday my good friend Bernie Maskit died.
Although not totally unexpected --  we all die --  a great shock nevertheless.  One more string cut. An end to our weekly Zoom meetings. He was my colleague (at MIT and Stony Brook and was my first or second appointment to the SB math department), collaborator (numerous joint papers), advisor (in mathematics and administration), but most of all friend (who else would drive me 60 miles, beginning at 3am, to pick up a delivery from abroad?). How did I ever reciprocate? He was a collaborative member of mathematics departments and a friend to many of his colleagues and students.
Hard to accept his loss, but nevertheless we must.

Maskit and I worked on similar problems mostly involving Kleinian groups and the Riemann surfaces associated with them, he from a topological and group theoretic approach and I as a complex analyst. Nevertheless we found a common language and rapidly developed joint programs, as is evident by the number of joint papers we produced. At the time we met, Kleinian groups and their moduli (deformations) were ready for a rebirth, mostly as a result of two important, independent, but related papers that had recently appeared, the first by Lipman Bers and the second by Bernard Maskit.  This revival included strong contributions (some in joint papers) by L.V. Ahlfors, L.Bers, C.J. Earle, I. Kra, A. Marden, B. Maskit, among others in the older generation.  More important, perhaps, was the number of young people attracted to this field.  Maskit was the authority on construction of complicated finitely generated Kleinian groups from simpler building blocks.  It was quite exciting for me to witness this act of creation and his interactions with his students and followers.
%\subsection*{Claude LeBrun:  Bernie Maskit remembered}
  \vskip5pt
\noindent\large{\bf Claude LeBrun:  Bernie Maskit remembered}
\vskip5pt
\normalsize
I was shocked and  deeply saddened when  my friend, collaborator, and colleague Bernie Maskit unexpectedly passed away last March.   Bernie had first befriended me in the early 1980s when I was an Assistant Professor at Stony Brook. He later became
chair of our department, and was serving in that capacity when I came up for tenure in 1988. 
Several  other Assistant Professors in the department had recently been denied tenure, so 
I was justifiably worried   about my own prospects. 
%the outcome of my own case. 
But Bernie's gentle, avuncular  encouragement
and   
generous  practical  advice helped   me somehow navigate and endure the  process without becoming a total nervous wreck, 
even as  he  helped shepherd my case  through  to a successful outcome. 

After I returned from my first sabbatical, Bernie and I started to meet for lunch from time to time, partly   to discuss mathematics, but also just to chat
about 
 life, the universe, and everything. My  work on the scalar curvature of self-dual $4$-manifolds had  convinced me of the importance 
of  locally-conformally-flat orbifolds,  but this  left me trying to puzzle out a host of subtle questions concerning Kleinian groups and their limit sets. 
 Fortunately, Bernie was a world-class expert on these matters, and he was only too happy to   educate
 me about  many of the  fine points I needed to understand.  It took  many years for  our discussions  to crystallize into a collaboration, but our joint
2008 paper,  {\it On Optimal $4$-Dimensional Metrics}, was well worth the wait. Our article  constructed 
Riemannian metrics that minimize the $L^2$-norm of the curvature tensor 
on connected sums of five or more copies of the complex projective
plane, and this then allowed us to prove 
a classification result that  I still find enormously satisfying. 

However,   soon after our paper was published,   Bernie  decided to retire and move into  Manhattan,  and I  saw him  only infrequently after that. 
But I still  miss his laughter, and 
I  still  miss our  conversations. 
 And I  am sure  that many of our  colleagues must also be  feeling the same sense of loss.  
%\subsection*{Dennis P. Sullivan}
  \vskip5pt
\noindent\large{\bf Dennis P. Sullivan}
\vskip5pt
\normalsize
I first met Bernie at MIT as a post doc and he  as a junior professor in 1970. I asked  him what he worked on for which he answered, ``I work on two by two matrices and surfaces."
 I hid my wonder since I was part of the field studying manifolds of dimension  five 
or more, which did not include three and four because they were too hard, and did not include two dimensions because they were too easy. (cf corrigendum below)

Later in the 70s, inspired by a question from Lipman Bers,  I learned about 
Kleinian groups and their limit sets and I tried to use the 
Klein-Maskit combination construction to answer Bers' question.
During this period I visited Stony Brook every summer and learned much more from Bernie and Irwin Kra about Schottky groups and automorphic forms, respectively.
I also learned how silly my first reaction at MIT was to Bernie's modest answer, ``I work on 2x2 matrices and surfaces,"  because 
the automorphism group of a surface of genus $g$ is a deep and rich mathematical object, actively studied, and central to several fields.

From Bernie I learned about the pleasant restaurants and cultural diversions of his hometown, Huntington, where he lived because his spouse worked in New York City  and he worked in   Stony Brook.

Bernie was smart, modest and likeable.
I am glad to have had the chance to know him. 

%\subsection*{Gadde  Swarup: RIP Bernie Maskit }
  \vskip5pt
\noindent\large{\bf Gadde  Swarup: RIP Bernie Maskit}
\vskip5pt
\normalsize
A very kind man who [helped me to] revive my research .... I had a brief collaboration with him and a longer collaboration was lost in the mail. I had many pleasant moments with him ... Mladen Bestvina tells me that I told him about the Maskit combination theorems, which led to his important work with Feighn on combination theorems for word hyperbolic groups.

\bibliographystyle{foo}
\bibliography{ExampleRefs}

\end{document}